\documentclass[a4paper]{amsart}

\usepackage{amssymb,amsmath,graphics}
\usepackage{latexsym}
\usepackage{eucal}
\makeatletter
\@addtoreset{equation}{section}\makeatother

\newtheorem{theo}{Theorem}[section]
\newtheorem{lem}[theo]{Lemma}
\newtheorem{prop}[theo]{Proposition}

\newcommand{\mc}{\mathcal}
\newcommand{\rr}{\mathbb{R}}
\newcommand{\nn}{\mathbb{N}}
\newcommand{\cc}{\mathbb{C}}
\newcommand{\hh}{\mathbb{H}}
\newcommand{\zz}{\mathbb{Z}}

\newcommand{\la}{\lambda}
\newcommand{\eps}{\epsilon}

\newcommand{\pl}{\partial}
\newcommand{\x}{\times}

\newcommand{\supp}{\textrm{supp}}
\newcommand{\cjd}{\rangle}
\newcommand{\cjg}{\langle}

\newcommand{\demi}{\frac{1}{2}}
\newcommand{\ndemi}{\frac{n}{2}}

\newcommand{\trans}{{^t}\!}

\def\qed{\hfill$\square$}
\begin{document}
\title[Resonances on some geometrically finite hyperbolic manifolds]
{Resonances on some geometrically finite hyperbolic manifolds}
\author[Colin Guillarmou]{Colin Guillarmou}
\address{Laboratoire de Math\'ematiques Jean Leray
         UMR 6629 CNRS/Universit\'e de Nantes (France)
         and Department of Mathematics, Purdue University, West Lafayette, IN (USA)}
     \email{cguillar@math.purdue.edu}
%
\begin{abstract}
\noindent We first prove the meromorphic extension of the resolvent to $\cc$ on
a class of geometrically finite hyperbolic manifolds with infinite
volume and we give a polynomial bound on the number
of resonances. This class notably contains the quotients $\Gamma\backslash\hh^{n+1}$
with rational non-maximal rank cusps previously studied by 
Froese-Hislop-Perry.
\end{abstract}
\maketitle

\section{Introduction}

The purpose of this work is to prove the meromorphic extension of
the resolvent as well as a polynomial bound of resonances for the
Laplacian on some geometrically finite hyperbolic manifolds
$\Gamma\backslash\hh^{n+1}$ whose non-maximal rank cusps are
`rational'. This condition is essentially equivalent to suppose
that the parabolic subgroups are conjugate to lattices of
translations acting on $\rr^n$.\\

Scattering theory, meromorphic continuation of the resolvent for
the Laplacian, Eisenstein series and distribution of resonances
have been deeply studied on geometrically finite hyperbolic
surfaces (see \cite{S,PA,CDV,G1,MU,GZ1,GZ3,Z1}). New
geometric difficulties appear in higher dimension, notably the
fact that a geometrically finite quotient
$\Gamma\backslash\hh^{n+1}$ is not, in general, a compact
perturbation of explicitly computable models and the
method used for surfaces can not be applied. However, when the manifold has a nice
structure near infinity, say when it conformally compactifies,
Mazzeo and Melrose \cite{MM} have found a powerful method to prove
the meromorphic continuation of the resolvent and to describe it
in details. Roughly, this conformal hypothesis is equivalent to
taking groups without parabolic elements. Perry \cite{P1,P2},
Joshi-Sa Barreto \cite{JSB}, Guillop\'e-Zworski \cite{GZ2}, and
Patterson-Perry \cite{PP} have studied the scattering matrix,
Eisenstein series and distribution of resonances for these classes
of manifolds and it is not difficult to see that parabolic 
elements with maximal rank can be added to the group without  
significative difficulties.   
Nevertheless, when the group has parabolic elements with non-maximal rank, 
some infinite volume cusps appear and most of those
cases remain quite mysterious in general, at least in the point of view of
the scattering theory and the meromorphic continuation of the
resolvent. It is worth noting that the analysis of the
spectrum of the Laplacian on forms and Hodge theory has been
written down by Mazzeo and Phillips \cite{MP} on  
geometrically finite hyperbolic manifolds and later, 
Froese-Hislop-Perry \cite{FHP1,FHP2} have 
studied scattering theory and Eisenstein functions on these manifolds 
in dimension $3$. 
Perry \cite{P} also proved the meromorphic continuation of the resolvent
in a small strip near the critical line for a class of manifolds
with rational non-maximal rank cusps in all dimensions.
Last but not least, Bunke and Olbrich \cite{BO} dealt with all cases 
of geometrically finite hyperbolic manifold using a very different approach, 
in particular they are able to extend the Eisenstein functions and scattering 
operator but they did not study the resolvent explicitly and have 
no result about the distribution of resonances.\\

In this work, we hope to give a rather simple way to prove the meromorphic 
continuation of the resolvent on a class of geometrically finite hyperbolic manifolds
which allows to bound the number of resonances (poles of the resolvent) 
in a disc of radius $R$ of $\cc$. 
Our method is very similar from the approaches of Froese-Hislop-Perry \cite{FHP2},
Perry \cite{P} or Guillop\'e-Zworski \cite{GZ2} in the sense
that it provides a precise analysis of the resolvent. 
Actually, the pseudo-differential structure of the resolvent and 
the scattering operator will be described in more 
details in \cite{G0}, the purpose of the present work
being essentially the estimate of the resonances distribution.\\

We thus consider an infinite volume hyperbolic manifold
$\Gamma\backslash \hh^{n+1}$ where $\Gamma$ is a discrete group of
isometries of $\hh^{n+1}$ which admits a fundamental domain with
finitely many sides (the manifold is said geometrically finite)
and such that each parabolic subgroup of $\Gamma$ is conjugate to
a lattice of translations in $\rr^n$. This is exactly the class of
manifolds studied by Perry \cite{P} and, as noticed in this work, 
it covers the case when no parabolic subgroup contains irrational
rotations, possibly by passing to a finite cover. The manifold
$X:=\Gamma\backslash \hh^{n+1}$ equipped with the hyperbolic
metric is then complete and the spectrum of the Laplacian
$\Delta_X$ splits into continuous spectrum
$[\frac{n^2}{4},\infty)$ and a finite number of eigenvalues
included in $(0,\frac{n^2}{4})$. Perry \cite{P} proved that the
modified resolvent
\[R(\la):=(\Delta_X-\la(n-\la))^{-1}\]
extends from $\{\Re(\la)>\ndemi\}$ to $\{\Re(\la)>\frac{n-1}{2}\}$
meromorphically with poles of finite multiplicity (i.e. the rank of
the polar part in the Laurent expansion at each pole is finite) in
weighted $L^2$ spaces.

We first show the
\begin{theo}\label{prolongmer}
Let $X=\Gamma\backslash \hh^{n+1}$ be a geometrically finite
hyperbolic manifold with infinite volume and such that each
parabolic subgroup of $\Gamma$ is conjugate to a lattice of
translations in $\rr^n$. Then the modified resolvent for the
Laplacian
\[R(\la)=(\Delta_X-\la(n-\la))^{-1}: L^2_{comp}(X)\to H^2_{loc}(X)\]
extends from $\{\Re(\la)>\ndemi\}$ to $\cc$ meromorphically with
poles of finite multiplicity.
\end{theo}

The poles of the resolvent are called resonances and are spectral
data which correspond in a sense to the eigenvalues of the compact
cases. They are closely related to the zeros and poles of
Selberg's zeta function. Our second result involves the asymptotic
distribution of resonances:

\begin{theo}\label{nbres}
With the  assumptions of Theorem \ref{prolongmer}, the number of
resonances $N(R)$ (counted with multiplicities) contained in the
complex disc $D(\ndemi,R)$ or radius $R$ satisfies
\[N(R)\leq CR^{n+2}+C\]
for some constant $C>0$.
\end{theo}

Notice that in a less general case, the non-optimal power
$n+2=\dim(X)+1$ already appeared in the work of Guillop\'e-Zworski
\cite{GZ2}. We also emphasize that our method extends to compact
perturbations of the hyperbolic manifolds considered in Theorems
\ref{prolongmer} and \ref{nbres}.

These theorems are proved by using a parametrix construction
of the resolvent. The main difficulty of the construction is 
that $X$ can not be splitted into one compact set and one neighbourhood
of infinity where a model resolvent is explicitly known. As for  
convex co-compact hyperbolic manifolds \cite{P1,P2,GZ2,PP}, we use a 
finite covering of a neighbourhood of the infinity $X(\infty)$ of $X$ 
by several model neighbourhoods with some exlicit formula 
for the resolvent of the laplacian on these models. 
Following the method of Perry \cite{P}, these model 
resolvents, cut-off by functions forming a partition of unity near infinity, 
give a first approximation of $R(\la)$ but it is not sufficient
to obtain a compact error. The important property we then 
use is that the error terms resemble those which appear 
in the convex co-compact manifolds (and more generally on 
asymptotically hyperbolic ones), in the sense that their Schwartz kernels
are smooth in $X\x X$ and have a polyhomogeneous asymptotic 
expansion at $X(\infty)\x X$ with support on the left factor which does 
not intersect the cusp points. Consequently the indicial equation used in \cite{MM,GZ2} allows to 
refine the parametrix of $R(\la)$ to get a residual term which is compact on all
$\rho^NL^2(X)$ ($N>0$) with $\rho$ a weight function converging to $0$ 
at infinity.\\

We will use the following notations: $\cjg
z\cjd:=(1+|z|^2)^{\demi}$; if $A$ is a compact operator on a
Hilbert space, $|A|:=(A^*A)^\demi$ and $(\mu_l(A))_l$ are the
eigenvalues of $|A|$ (called the singular values of $A$); $C$ will
denote a large constant, not necessarily always the same.
We shall also often identify operators with their Schwartz kernels.\\

\textbf{Acknowledgement.} This work has been begun at Nantes University
and finished at Purdue University. I would like to thank Peter Perry, Rafe Mazzeo 
and Martin Olbrich for pointing out to me some references about the subject.

\section{Geometry of the manifold}

This section is strongly inspired by Perry's paper \cite{P}, the reader can refer 
to it for more details (see also  \cite{MP,FHP1,FHP2}).
We consider a hyperbolic quotient $X=\Gamma\backslash \hh^{n+1}$ where
$\Gamma$ is a geometrically finite discrete group of hyperbolic isometries
with no elliptic elements, and such that all parabolic subgroups are
conjugate to lattices of translations acting on $\rr^n$.
In this case there exists a compact $K$ of $X$ such that $X\setminus K$
is covered by a finite number of charts isometric to either
a regular neighbourhood $(M_r,g_r)$ or a rank-k cusp neighbourhood
$(M_k,g_k)$ (with $1\leq k\leq n$), where
\[M_r=\{(x,y)\in(0,\infty)\x\rr^n; x^2+|y|^2<1\}, \quad g_r=x^{-2}(dx^2+dy^2),\]
\[M_k=\{(x,y,z)\in (0,\infty)\x\rr^{n-k}\x T^k; x^2+|y|^2>1\},\quad
g_k=x^{-2}(dx^2+dy^2+dz^2)\]
for $k<n$ with $(T^k,dz^2)$ a k-dimensional flat torus and
\[M_n=\{(x,z)\in (0,\infty)\x T^n; x>1\},\quad
g_n=x^{-2}(dx^2+dz^2)\]
with $(T^n,dz^2)$ a n-dimensional flat
torus. For notational simplicity, we will make as if there was one
neighbourhood of each type. There exist some smooth functions
$\chi^i,\chi^r,\chi^1,\dots,\chi^n$ on respectively $X,M_r,M_1,\dots,M_n$
which, through
the isometric charts $I_r,I_1,\dots,I_n$, satisfy
\begin{equation}\label{partition}
I_r^*\chi^r+\sum_{k=1}^nI_k^*\chi^k+\chi^i=1
\end{equation}
with $\chi^i$ having compact support in $X$.

We will also use cutoff functions in what
follows, thus we define
\begin{equation}\label{phi1}
\phi,\phi_L\in
C_0^\infty([0,2)), \quad \phi_L=1 \textrm{ on }[0,1],\quad
\phi=1 \textrm{ on } \supp(\phi_L). \end{equation}

\subsection{The non-maximal rank cusps neighbourhoods}

Let $X_k=\Gamma_k\backslash\hh^{n+1}$ be the quotient of
$\hh^{n+1}$ by a rank-k parabolic subgroup $\Gamma_k$ of $\Gamma$ which fixes
a single point at infinity of $\hh^{n+1}$. Modulo conjugation by a
 hyperbolic isometry, one can suppose that the fixed point is the point at 
infinity of $\hh^{n+1}$ in the half-space model $(0,\infty)\x\rr^n$.
$\Gamma_k$ can then be considered as a rank-k lattice of translations acting on 
$\rr^n$ (actually on a subspace of $\rr^n$ isomorphic to $\rr^k$), 
therefore it is the image of the lattice $\zz^k$ by a map $A_k\in GL_k(\rr)$ and the
flat torus $T^k:=\Gamma_k\backslash \rr^k$ is well defined.
$X_k$ is isometric to $\rr^+_x\x \rr_y^{n-k}\x T^k_z$ equipped
with the metric
\[g_k=\frac{dx^2+dy^2+dz^2}{x^2}\]
$dz^2$ being the flat metric on a k-dimensional torus $T^k$.
$M_k$ is the subset of $X_k$ with $x^2+|y|^2>1$
and we will often consider $\rr^+\x \rr^{n-k}$ as the $n-k+1$-dimensional
hyperbolic space $\hh^{n-k+1}$. Hence,
\begin{equation}\label{rhok}
\rho_k(x,y,z)=\rho_k(x,y):=\frac{x}{|y|^2+x^2+1}=
(2\cosh(d_{\hh^{n-k+1}}(x,y;1,0)))^{-1}\end{equation}
will be a natural weight function on $M_k$.

\begin{figure}[ht!]
\begin{center}
\input{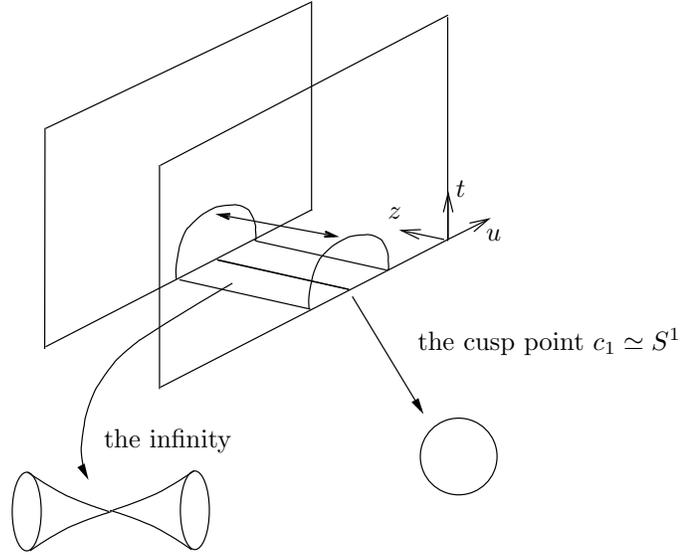}
\caption{The compactified neighbourhood of the cusp point $c_1$ in dimension $3$}
\end{center}
\end{figure}

Following Perry \cite{P}, we remark that the change of coordinates
\begin{equation}\label{coordtu}
t:=\frac{x}{x^2+|y|^2}, \quad u:=\frac{-y}{x^2+|y|^2}
\end{equation}
is an isometry from $M_k$ to
\[\{(t,u,z)\in(0,\infty)\x\rr^{n-k}\x T^k; t^2+|u|^2<1\}\]
equipped with the metric
\begin{equation}\label{metriccompact}
\frac{dt^2+du^2+(t^2+|u|^2)^2dz^2}{t^2}
\end{equation}
and $\rho_k(t,u)=\rho_k(x,y)$. 
These coordinates compactify $M_k$ and the infinity of $X$ in this
neighbourhood is then given by $\{\rho_k=0\}$ or equivalently
$\{t=0\}$. The `cusp point' here becomes a torus $c_k:=\{t=u=0\}\simeq T^k$. 

Without loss of generality and possibly by adding
regular charts in the covering of a neighbourhood of the infinity of $X$,
we can choose the cut-off function for the rank-k cusp neighbourhood such that
$\chi^k(x,y):=1-\phi(x)\psi^k(y)$ with $\psi^k(y)\in
C_0^\infty(|y|_{\rr^{n-k}}<2)$ and $\psi^k(y)=1$ on
$\{|y|\leq 1\}$. We also set $\psi^k_L\in C_0^\infty(|y|<2)$ such
that $\psi^k=1$ on $\supp(\psi_L^k)$,
$\psi_L^k=1$ on $\{|y|\leq 1\}$ and we define 
$\chi^k_L(x,y):=1-\phi_L(x)\psi^k_L(y)$, which satisfies 
$\chi_L^k=1$ on $\supp(\chi^k)$.

\subsection{The maximal rank cusps}
Let $X_n=\Gamma_n\backslash\hh^{n+1}$ be the quotient of $\hh^{n+1}$
by a rank-n parabolic group subgroup $\Gamma_n$ of $\Gamma$, it is then isometric to
$\rr^+_x\x T^n_z$ equipped with the metric
\[g_n=\frac{dx^2+dz^2}{x^2}\]
$dz^2$ being the flat metric on the n-dimensional torus $T^n=\Gamma_n\backslash\rr^n$.
As before $M_n$ is the subset of $X_n$ with $x>1$ and the weight function we choose on
$M_n$ is $\rho_n=x^{-1}$.
Taking $u=x^{-1}$, $M_n$ is also isometric to
\[M_n=\{(u,z)\in (0,\infty)\x T^n; u<1\}\]
with the metric
\[u^{-2}(du^2+u^4dz^2)\]
and $\rho_n=u$. The infinity of $X$ in this neighbourhood
is given by the `cusp point' $c_n:=\{\rho_n=0\}\simeq T^n$ which is a torus.

Here $\chi^n$ can be taken depending
only on $x$, for example $\chi^n(x,z)=1-\phi(x)$ and we set
$\chi_L^n:=1-\phi_L(x)$, hence $\chi_L^n=1$ on $\supp(\chi^n)$.

\subsection{The regular neighbourhoods}
The regular neighbourhoods
are those treated in \cite{GZ2} and will not
be discussed in details.
We just recall that the weight function here is
$\rho_r:=x$ and that the infinity of $X$ in this neighbourhood
is given by $\rho_r=0$.

The function $\chi^r$ can be chosen so that
$\chi^r=\phi_L(\frac{x}{\eps})\psi^r(y)$ with $\psi^r\in C_0^\infty(|y|<1)$ and $\eps>0$ small.
Let $\psi_L^r\in C_0^\infty(|y|<1)$
with $\psi_L^r=1$ on $\supp(\psi^r)$
and $\chi_L^r(x,y)=\phi_L(\frac{x}{2\eps})\psi_L^r(y)$, we thus have
$\chi_L^r=1$ on $\supp(\chi^r)$.

\subsection{Weight function, compactification}
We can then define as weight function the following
\begin{equation}\label{weight}
\rho:=\chi^i+I_r^*(\chi^r\rho_r)+\sum_{k=1}^nI_k^*(\chi^k\rho_k).
\end{equation}
$X$ can be compactified in a compact manifold with boundary
$\bar{X}$ such that $\rho$ is a boundary defining function of $\bar{X}$.
The boundary can be decomposed in the form
$\pl\bar{X}=\bar{B}\sqcup T^n$ with $\bar{B}$ a smooth compact manifold
and $T^n$ the torus coming from the maximal rank cusp.
From the discussion above, we see that the metric on $X$ can
be expressed by
\[g=\frac{H}{\rho^2}\]
with $H$ a smooth non-negative symmetric 2-tensor on
$\bar{X}$ which degenerates at cusps points $(c_k)_{k=1,\dots,n}$.
Let $B:=\bar{B}\setminus \{c_1,\dots,c_n\}$, the restriction $H|_{B}$
is then a smooth metric on the non-compact manifold $B$.
In a sense $B$ will be the geometric infinity where the scattering
can occur.

\begin{figure}[ht!]
\begin{center}
\input{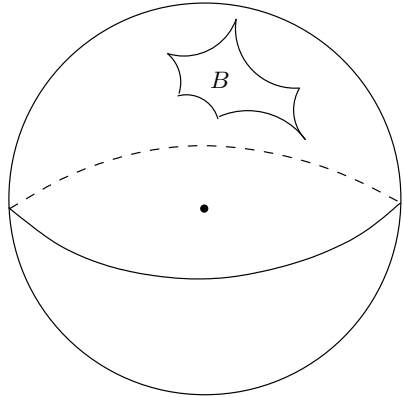}
\caption{A fundamental domain in $S^n$ for the infinity $B$ of $X=\Gamma\backslash \hh^{n+1}$}
\end{center}
\end{figure}

We need to construct a new weight function depending on a small
paramater $\delta>0$, which is equal to $1$ on a big compact (depending
on $\delta$) and to $\rho$ near infinity.

For $\alpha=1,\dots,n,r$,  we have $\rho_\alpha(w)={I_\alpha}_*\rho(w)$
when $\chi^\alpha(w)=1$ and there exists $C>1$
such that
\begin{equation}\label{comprhok}
C^{-1}\rho_\alpha\leq {I_\alpha}_*\rho\leq C\rho_\alpha \textrm{ on }
\supp(\chi^\alpha).
\end{equation}
Let $1>\delta>0$ be small, we define the new weight function on $X$
\[\rho_\delta(w):=(1-\phi_L(4\rho(w)\delta^{-1}))+
\rho(w)\phi_L(4\rho(w)\delta^{-1})\]
and using that $\rho\leq 1$ we can check that
 $\rho_\delta=\rho$ on
$\{\rho\leq \frac{\delta}{4}\}$, $\rho_\delta=1$ on
$\{\rho>\frac{\delta}{2}\}$ and $\rho_\delta\geq \rho$ everywhere,
hence
\begin{equation}\label{comprho}
1\leq \frac{\rho_\delta}{\rho}\leq
\frac{4}{\delta}.
\end{equation}

We also define the function
\[\chi^\alpha_\delta(w):=\phi_L(\rho(w)\delta^{-1})
\chi^\alpha(w).\]
Thus $\chi^\alpha_L=1$ on the support of $\chi^\alpha_\delta$ and
\begin{equation}\label{rhomajdelta}
{I_\alpha}_*\rho\leq 2\delta
\textrm{ on } \supp(\chi^\alpha_\delta).
\end{equation}
In view of (\ref{partition}), we deduce
\[I_r^*\chi_\delta^r+\sum_{k=1}^nI_k^*\chi_\delta^k+\chi^i_\delta=1 \]
with
\[\chi^i_\delta(w):=1-\phi_L(\rho(w)\delta^{-1})+\chi^i(w)
\phi_L(\rho(w)\delta^{-1})\]
having a compact support included in $\{\rho\geq \delta\}$
for $\delta$ sufficiently small.
Finally let
\[\chi^i_{L,\delta}(w):=1-\phi_L(2\rho(w)\delta^{-1})\]
which is a smooth function with support in
$\{\rho\geq \frac{\delta}{2}\}$ and
$\chi^i_{L,\delta}=1$ on $\supp(\chi^i_\delta)$, thus $\rho_\delta=1$ on
$\supp(\chi^i_{L,\delta})$.

\begin{figure}[ht!]
\begin{center}
\input{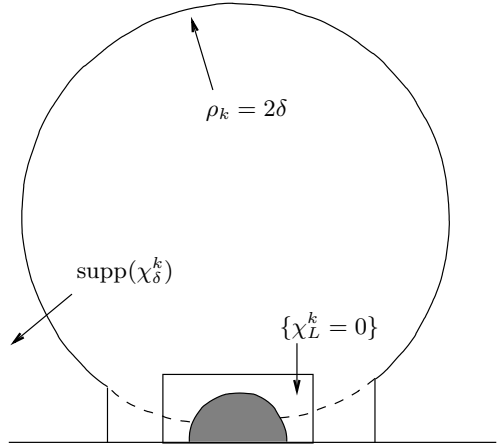}
\caption{The support of cut-off functions in $M_k$}
\end{center}
\end{figure}

Later, this parameter $\delta$ will be chosen small
enough to insure that our residual terms in the parametrix
construction have a small norm for $\Re(\la)\gg \ndemi$.
This is, in a sense, the idea used in \cite{GZ2}.
For what follows and for simplicity, 
we will often write $\rho_\delta$ instead of ${I_k}_*\rho_\delta$.

\section{Parametrix and estimates}\label{parametestim}
\subsection{The non-maximal rank cusps}
We recall that $X_k=\Gamma_k\backslash \hh^{n+1}$ with $\Gamma_k$
the image of the lattice $\zz^k$ by a map $A_k\in GL_k(\rr)$. By
using a Fourier decomposition on the torus
$T^k=\Gamma_k\backslash\rr^k$ and conjugating by
$x^{\frac{k}{2}}$, the operator $\Delta_{X_k}-\la(n-\la)$ acts on
\[L^2(X_k)=\bigoplus_{m\in\zz^k}\mc{H}_m,\quad
\mc{H}_m\simeq L^2(\rr_y^{n-k}\x\rr^+_x,x^{-(n-k+1)}dydx)=L^2(\hh^{n-k+1})\]
as a family of operators
\[P_m(\la):=-x^2\pl_x^2+(n-k-1)x\pl_x+x^2(\Delta_y+|\omega_m|^2)-s(n-k-s)\]
where $\omega_m=2\pi \trans (A_k^{-1})m$ for $m\in\zz^k$ are the eigenvalues
of the Laplacian on $T_k$ (with eigenfunction $e^{i\omega_m.z}$)
and $s:=\la-\frac{k}{2}$ is a shifted spectral parameter.
The resolvent $R_{X_k}(\la)=(\Delta_{X_k}-\la(n-\la))^{-1}$
for the Laplacian on $X_k$ is computed by Perry \cite{P}
for $\Re(\la)>\ndemi$
\begin{equation}\label{resolventxk}
R_{X_k}(\la)=\bigoplus_{m\in\zz^k}R_m(\la) \textrm{ on }
L^2(X_k)=\bigoplus_{m\in\zz^k}\mc{H}_m
\end{equation}
with
\begin{equation}\label{fourier}
R_m(\la;x,y,x',y')=|A_k|^{-\demi}(xx')^{-\frac{k}{2}}
\int_{\rr^k}R_{\hh^{n+1}}(\la;x,y,z;x',y',0)e^{i\omega_m.z}dz
\end{equation}
where $R_{\hh^{n+1}}(\la)=(\Delta_{\hh^{n+1}}-\la(n-\la))^{-1}$ is
the resolvent for the Laplacian on $\hh^{n+1}$ and
$|A_k|:=|\det(A_k)|$. We set
\begin{equation}\label{rdtau}
 r=(|y-y'|^2+x^2+{x'}^2)^\demi,\quad d=\frac{xx'}{r^2},\quad
\tau=\frac{xx'}{r^2+|z|^2}=d(1+\frac{|z|^2}{r^2})^{-1}
\end{equation}
and recall (see e.g. \cite{GZ2}, \cite{P}) that the resolvent on $\hh^{n+1}$ can be written for all $J\in\nn\cup\infty$ 
\begin{equation}\label{reshyp}
R_{\hh^{n+1}}(\la;x,y,z;x',y',0)=\tau^\la\sum_{j=0}^{J-1}\alpha_{j,n}(\la)\tau^{2j}
+\tau^{\la+2J}G_{J,n}(\la,\tau)
\end{equation}
\[\alpha_{j,n}(\la):=\frac{2^{-1}\pi^{-\ndemi}\Gamma(\la+2j)}
{\Gamma(\la-\ndemi+1+j)\Gamma(j+1)}\]
with $G_{J,n}(\la,\tau)$ a smooth function in $\tau\in[0,\demi)$ with a conormal 
singularity at $\tau=\demi$ and $G_{\infty,n}(\la,\tau)=0$.
Note that the sum (\ref{reshyp}) converges locally uniformly in $\tau\in[0,\demi)$
if $J=\infty$.
From (\ref{fourier}) and (\ref{reshyp}) it is easy to see, by the change of
variable $w=z/r$, that for $m\not=0$ and setting $s:=\la-\frac{k}{2}$ 
\begin{equation}\label{rkm}
R_m(\la)=d^{s}\sum_{j=0}^{J-1} \alpha_{j,n}(\la)d^{2j}F_{j,\la}(
r\omega_m)+d^{s+2J}\int_{\rr^k}e^{-ir\omega_m.z}\frac{G_{J,n}(\la,d(1+|z|^2)^{-1})}{(1+|z|^2)^{\la+2J}}dz
\end{equation}
\begin{equation}\label{fjla}
F_{j,\la}(u)=|A_k|^{-\demi}\int_{\rr^k}e^{-iu.w}(1+|w|^2)^{-\la-2j}dw=
|A_k|^{-\demi}\frac{2^{-\la-2j+1}(2\pi)^{\frac{k}{2}}}{\Gamma(\la+2j)}|u|^{s+2j}
K_{-s-2j}(|u|)
\end{equation}
when $\Re(\la)>\ndemi$ (see e.g. \cite{GS} for the last formula), 
$K_s(z)$ being the Bessel function defined by
\begin{equation}\label{besselk}
K_s(z):=\int_0^\infty\cosh(st)e^{-z\cosh(t)}dt, \quad z>0.
\end{equation}
It is easy to see (and will be studied later) 
that the sum (\ref{rkm}) with $J=\infty$ converges uniformly for $r>0$ and $d\in[0,\demi)$.
When $m=0$, $R_0(\la)$ is the shifted Green kernel of the
Laplacian on $\hh^{n-k+1}$, that is
\begin{equation}\label{rk0}
R_0(\la)=d^s\sum_{j=0}^{J-1} \alpha_{j,n-k}(s)d^{2j}+d^{s+2J}G_{J,n-k}(\la,d), \quad
s=\la-\frac{k}{2}.
\end{equation}
For simplicity, we will write (\ref{rk0})
under the form (\ref{rkm}) by defining
\[ F_{j,\la}(0):=\frac{\alpha_{j,n-k}(s)}{\alpha_{j,n}(\la)}.\]
The representations (\ref{rkm}) and (\ref{rk0}) give a
meromorphic extension of $R_m(\la)$ to $\cc$, with poles on
$\frac{k}{2}-\nn_0$ of finite rank which only come from the case $m=0$
when $n-k+1$ is even. The continuity property of the extented 
operators on weighted $L^2$ spaces will be checked later.

We are now going to find a parametrix for $\Delta_X-\la(n-\la)$ on
the neighbourhood $I_k^{-1}(M_k)$ of our manifold $X$, by using
the explicit formulae given before for the cusp $X_k$. The
constructions are very similar to those of Guillop\'e-Zworski
\cite{GZ2} for the conformally compact ends. After the Fourier
decomposition, the construction of a parametrix will be obtained
on each $\mc{H}_m$ from the model resolvent and an iterative
process as in \cite{GZ2} but for a $n-k+1$ dimensional hyperbolic
Laplacian with the potential $x^2|\omega_m|^2$.

\begin{prop}\label{parametrixkcusp}
For $N\in\nn$ large, there exist some bounded operators
\[\mc{E}^k_{N}(\la): \rho_\delta^NL^2(M_k)\to \rho_\delta^{-N}L^2(M_k)\]
\[\mc{K}^k_N(\la): \rho_\delta^NL^2(M_k)\to \rho_\delta^NL^2(M_k)\]
meromorphic in $\{\Re(\la)>-\frac{N+1}{2}\}$ with simple poles
at $\frac{k}{2}-j$ (with $j\in\nn_0$) of ranks uniformly
bounded by $C(j+1)^{n-k+1}$ such that
\[(\Delta_{M_k}-\la(n-\la))\mc{E}^k_N(\la)=\chi_\delta^k+\mc{K}^k_N(\la),\]
$\mc{K}^k_N(\la)$ is trace class on $\rho_\delta^NL^2(M_k)$ and for
$|\la|\leq \frac{N}{2}$, $\textrm{dist}(\la,\frac{k}{2}-\nn_0)>\frac{1}{8}$
and $q>0$ there exists $C_{\delta,q}>0$ such that
\begin{equation}\label{detkcusp}
\det(1+q|\mc{K}^k_N(\la)|)\leq e^{C_{\delta,q}\cjg N\cjd^{n+2}},
\end{equation}
the determinant being on $\rho_\delta^NL^2(M_k)$.
Moreover for $\la_N=\frac{N}{4}$,
\begin{equation}\label{normlan}
||\mc{K}^k_N(\la_N)||_{\mc{L}(\rho_\delta^NL^2(M_k))}\leq
(C\delta)^{\frac{N}{4}}. \end{equation}
\end{prop}
\textsl{Proof}: we first set
\[E_0(\la):=\chi_L^kR_{X_k}(\la)\chi^k_\delta,\quad
K_0(\la):=\phi_L[\Delta_{X_k},\psi^k_L]R_{X_k}(\la)\chi^k_\delta,\quad
L^\sharp(\la):=\psi_L^k[\Delta_{X_k},\phi_L]R_{X_k}(\la)\chi^k_\delta\]
and, since $\Delta_{M_k}=\Delta_{X_k}$ as a differential operator on $M_k$, we clearly have
\[(\Delta_{M_k}-\la(n-\la))E_0(\la)=\chi^k_\delta+K_0(\la)+L^\sharp(\la).\]
Since the functions in the range of $L^\sharp(\la)$ have compact support,
$L^\sharp(\la)$ is compact on our weighted spaces but $K_0(\la)$ is not.
However, it is important to note that the range of $K_0(\la)$ is composed
of functions whose support does not intersect the cusp point, thus they
can be included in a regular neighbourhood of infinity and the iterative method of
Mazzeo-Melrose \cite{MM} (or \cite{GZ2}) can then be used to remove all the
Taylor expansion of $K_0(\la;w,w')$ at the boundary $\{\rho(w)=0\}$.

To achieve it, one can decompose $E_0(\la)=\oplus_{m\in\zz^k}E_{0,m}(\la)$ and
$K_0(\la),L^\sharp(\la)$ similarly
on the Fourier modes of $T^k$ and using the new variables
$v:=x^2$ and $v':={x'}^2$ as in \cite{GZ2} we can write
\[K_{0,m}(\la)=\phi_L\sum_{j=0}^{2N}K^j_{0,m}(\la)+K^{\sharp}_{2N,m}(\la)\]
where for $j=0,\dots,2N$, $w=(v,y)$, $w'=(v',y')$, $r=r(w,w')$ and $d=d(w,w')$
\[K^j_{0,m}(\la;w,w')=v^{\frac{s}{2}+j+1}[\Delta_y,\psi_L^k]r^{-2s-4j}
\alpha_{j,n}(\la)F_{j,\la}(r\omega_m)\chi^k_\delta(w'){v'}^{\frac{s}{2}+j}\]
\[K^\sharp_{2N,m}(\la;w,w')=
v\phi_L[\Delta_y,\psi_L^k]d^s\sum_{j=2N+1}^\infty\alpha_{j,n}(\la)d^{2j}F_{j,\la}(
r\omega_m)\chi^k_\delta(w')\]
At last we complete the construction by mimicking \cite[Prop. 3.1]{GZ2}
\begin{equation}\label{defm}
E_j(\la)=\bigoplus_{m\in\zz^k}E_{j,m}(\la), \quad
K_j(\la)=\bigoplus_{m\in\zz^k}K_{j,m}(\la),\quad
L_j(\la)=\bigoplus_{m\in\zz^k}L_{j,m}(\la)
\end{equation}
with the induction formulae for $j=1,\dots,2N$
\[L_{0,m}(\la)=0,\]
\[E_{j,m}(\la):=E_{j-1,m}(\la)+[2j(2\la+2j-n)]^{-1}K^{j-1}_{j-1,m}(\la),\]
\[K_{j,m}(\la):=K^j_{j,m}(\la)+K^{j+1}_{j,m}(\la)+\dots+K^{2N}_{j,m},\]
\[K^j_{j,m}(\la):=K^j_{j-1,m}(\la)+[2j(2\la+2j-n)]^{-1}Q_m(s/2+j)K^{j-1}_{j-1,m}(\la),\]
\[K^l_{j,m}(\la):=K^l_{j-1,m}(\la) \textrm{ for } l=j+1,\dots,2N\]
\[L_{j,m}(\la):=L_{j-1,m}(\la)+[2j(2\la+2j-n)]^{-1}[\Delta_{X_k},\phi_L]
K^{j-1}_{j-1,m}(\la),\]
where $Q_m(\zeta)$ is defined by
\[Q_m(\zeta):=2(n-k-2-4\zeta)\pl_v+4v\pl^2_v+\Delta_y+|\omega_m|^2.\]
Using the crucial 'indicial relation' (see \cite[Eq. 3.12]{GZ2})
\[(\Delta_{\hh^{n-k+1}}+v|\omega_m|^2)v^{\zeta}f(v,y)=
2\zeta(n-k-2\zeta)v^{\zeta}f(v,y)+v^{\zeta+1}Q_m(\zeta)f(v,y),\]
we then obtain from the previous construction that
\[P_m(\la)E_{2N,m}(\la)=\chi^k_\delta+\phi_LK_{2N,m}(\la)+L_{2N,m}(\la)
+K^\sharp_{2N,m}(\la)+L^\sharp_m(\la)\]
and a straightforward calculus shows that for $w=(v,y),w'=(v',y')$
\[K_{2N,m}(\la;w;w')=v^{\frac{s}{2}+2N+1}\sum_{j=0}^{2N-1}\beta_{j,2N}(\la)
\prod_{k=1}^{2N-j}Q_m(\frac{s}{2}+j+k)_{w}H_{j,m}(\la;w,w'),\]
\[L_{2N,m}(\la;w,w')=[\Delta_{X_k},\phi_L]
\sum_{0\leq j\leq p\leq 2N}\beta_{j,p}(\la)v^{\frac{s}{2}+p}
\prod_{k=1}^{p-j}Q_m(\frac{s}{2}+j+k)_wH_{j,m}(\la;w,w')\]
with
\[\beta_{j,p}(\la):=\frac{\pi^{-\ndemi}\Gamma(\la+2j)2^{-1}}{\Gamma(\la-\ndemi+1+p)
\Gamma(p+1)}\]
\[H_{j,m}(\la;w,w')=[\Delta_y,\psi_L^k]r^{-2s-4j}
F_{j,\la}(r\omega_m){v'}^{\frac{s}{2}+j}\chi^k_\delta(w').\]
The distributional kernels of
$E_{2N,m}(\la),K_{2N,m}(\la),L_{2N,m}(\la)$ and
$K^\sharp_{2N,m}(\la)$ are holomorphic in $\cc$ when $m\not=0$ and
meromorphic with simple poles of finite rank at each
$\frac{k}{2}-j$ ($j\in\nn_0$) when
$\omega_m=0$, the ranks being bounded by $C(1+j)^{n-k+1}$
(see again \cite[Prop. 3.1]{GZ2} for details).\\

At this stage we can set
\[\mc{E}^k_N(\la):=\bigoplus_{m\in\zz^k}E_{2N,m}(\la),\]
\[\mc{K}^k_N(\la):=\bigoplus_{m\in\zz^k}\Big(\phi_LK_{2N,m}(\la)+L_{2N,m}(\la)
+K^\sharp_{2N,m}(\la)+L^\sharp_m(\la)\Big)=\bigoplus_{m\in\zz}\mc{K}^k_{N,m}(\la).\]
Notice that the sums (\ref{defm}) are just formal so far, but we will
show their convergence in the following lemmas.

We will first show that $\mc{K}^k_{N,m}(\la)$ is compact on
$\rho_\delta^N\mc{H}_m$ and we will bound its singular values uniformly with
respect to $m$.
These estimates will prove that $\bigoplus_{m=0}^M\mc{K}^k_{N,m}(\la)$
converges to an operator $\mc{K}^k_{N}(\la)$ compact on
$\rho_\delta^NL^2(X_k)$ when $M\to\infty$, whose singular values are
$(\mu_l(\mc{K}^k_{N,m}(\la)))_{l\in\nn,m\in\zz^k}$ if
$(\mu_l(\mc{K}^k_{N,m}(\la)))_{l\in\nn}$ are the
singular values of $\mc{K}^k_{N,m}(\la)$ on $\rho_\delta^N\mc{H}_m$.

\begin{lem}\label{singvalk}
The operators $\phi_LK_{2N,m}(\la)$ and $L_{2N,m}(\la)$ are
trace class on $\rho_\delta^N\mc{H}_m$ for $|\la|\leq
\frac{N}{2}$, $\textrm{dist}(\la,\frac{k}{2}-\nn_0)>\frac{1}{8}$ and
their singular values satisfy
\[\mu_l\Big(\phi_LK_{2N,m}(\la)+L_{2N,m}(\la)\Big)\leq
e^{-\frac{\eps_0\cjg\omega_m\cjd}{4}}
\delta^{\Re(s)}\left(Cl^{-\frac{1}{n-k+1}}N\right)^{2N}
\max\Big(1,\big(\frac{\cjg\omega_m\cjd}{N}\big)^{\Re(s)-|\Re(s)|}\Big)\]
there for some $\eps_0>0$. Moreover, for $\la_N=\frac{N}{4}$, we have
\begin{equation}\label{majklnorm}
\left|\left|\phi_LK_{2N}(\la_N)+L_{2N}(\la_N)\right|\right|_{\mc{L}(\rho^N_\delta L^2(M_k))}\leq
(C\delta)^{\frac{N}{4}}
\end{equation}
\end{lem}
\textsl{Proof}: to begin, we give for $\eps>0$ an estimate
on the Bessel function in $\Re(z)>2\eps$
\begin{equation}\label{ks}
e^{\Re(z)}|K_s(z)|\leq C_\eps\sup_{t\geq 0}
[e^{-(\Re(z)-\eps)(e^t-1)+|\Re(s)|t)}]\leq
C_\eps\max\left(1,\left(\frac{\cjg\Re(s)\cjd}{|\Re(z)-\eps|}\right)^{|\Re(s)|}\right).
\end{equation}
To see that $\phi_LK_{2N,m}(\la)$ is trace class on $\rho_k^N\mc{H}_m\simeq\rho_\delta^N\mc{H}_m$,
we use a standard trick. Let $\Omega\subset \rr^{n-k+1}$ be an open ball
containing $\{x^2+|y|^2\leq 4\}$ and $\Delta_\Omega$ the Dirichlet
realization of the Laplacian on $\Omega$. Since
$(\Delta_\Omega+1)^{-M}$ is trace class for $M>\frac{n-k+1}{2}$ on
$L^2(\Omega)$, it suffices to show that
$(\Delta_{\Omega}+1)^N\rho_k^{-N}\phi_LK_{2N,m}(\la)\rho_k^N$ can be
extended as a bounded operator on $L^2(\Omega)$, and a uniform
bound on its norm together with a comparison with the singular
values of $(\Delta_\Omega+1)^{-N}$ will give an estimate for the
singular values $(\mu_l(\phi_LK_{2N,m}(\la)))_{l\in\nn,m\in\zz^k}$ on $\rho_k^N\mc{H}_m$. The
same method can be applied for $L_{2N,m}(\la)$.

By Stirling's formula and the complement formula, we check that for
$p\geq 0$ and $|\la|<\frac{N}{2}$
\begin{equation}\label{estimalphaj}
\left|\beta_{j,p}(\la)\frac{2^{-\la-2j+1}(2\pi)^{\ndemi}}{\Gamma(\la+2j)}\right|
\leq C^{N+p}p^{-p}\cjg p+\la\cjd^{-p-\Re(\la)}\leq
C^{N+p}N^{-2p-\Re(\la)}. \end{equation}
A straightforward estimate for
$|\la|\leq \frac{N}{2}, |\alpha|\leq 2N$ shows that
\[|\pl^\alpha_w(v^{\frac{s}{2}+2N+1}\rho_k(v,y)^{-N})|\leq
C^N(|\alpha|+N)^{|\alpha|}\] for $w=(v,y)\in\{0\leq v\leq 2,
|y|\leq 2\}$. We now choose the cut-off functions
$\psi_L^k,\phi_L$ quasi-analytic of order $5N$, that is
\begin{equation}\label{quasi}
||\pl_y^\alpha\psi_L^k(y)||_\infty\leq
(CN)^{|\alpha|}, \quad ||\pl_x^l\phi_L(x)||_\infty \leq
(CN)^{l}\quad \textrm{ for } |\alpha|\leq 5N, l\leq 5N.
\end{equation}
Therefore, for all smooth function $f(v,y)$ with support in
$\{v\in[0,2],|y|\leq 2\}$, we have for $M\leq N$ and $p\leq 2N$
\begin{equation}\label{estimprodq}
\left|(\Delta_w+1)^M\phi_L(w)\rho_k(w)^{-N}v^{\frac{s}{2}+2N+1}
\prod_{k=1}^{2N-j}Q_m(\frac{s}{2}+j+k)f(v,y)\right|_\infty
\end{equation}
\[\leq C^N\sum_{l_0+l_1+|l_2|+|l_3|=2N+M-j}\cjg\omega_m\cjd^{2l_0}(CN)^{2l_1+|l_3|}
|\pl^{2l_2+l_3}_wf(w)|_\infty\] and
\begin{equation}\label{estimprodq2}
\left|(\Delta_w+1)^M[\Delta_{X_k},\phi_L]\rho_k(w)^{-N}v^{\frac{s}{2}+p+1}
\prod_{k=1}^{p-j}Q_m(\frac{s}{2}+j+k)_{w}f(v,y)\right|_\infty
\end{equation}
\[\leq C^N\sum_{l_0+l_1+|l_2|+|l_3|=p+M-j}\cjg\omega_m\cjd^{2l_0}(CN)^{2l_1+|l_3|}
|\pl^{2l_2+l_3}_wf(w)|_\infty\]
(recall that $[\Delta_{X_k},\phi_L]$ has compact support). For $w'$ fixed,
we want to extend $r(.,w')$ in a complex neighbourhood of $\rr^{n-k+1}$ in $\cc^{n-k+1}$
to obtain bounds on its derivatives by Cauchy formula.
Using the fact that $r(w,w')>\eps_0$ for some
$\eps_0>0$ when $w\in\textrm{supp}\nabla\chi^k_L$,
$w'\in\textrm{supp}(\chi^k_\delta)$ we argue that for $(v_0,y_0)\in\rr\x\rr^n$ such that
$|v_0|+|y_0|\leq \eps<\demi\eps_0$ then
\[r^2=r(v+iv_0,y+iy_0,w')^2=(v+v'+|y-y'|^2-|y_0|^2)+i(v_0+2(y-y').y_0)\]
satisfies in $w=(v,y)\in\textrm{supp}\nabla\chi^k_L$,
$w'=(v',y')\in\supp(\chi^k_\delta)$
\begin{equation}\label{rcomplex}
\Re(r^2)>\demi\eps_0^2, \quad \textrm{arg}(r^2)\leq C\eps.
\end{equation}
This implies, with (\ref{ks}), that for $|\la|\leq \frac{N}{2}$,
$\textrm{dist}(\la,\frac{k}{2}-\nn_0)>\frac{1}{8}$, the function
\[\theta_{j,\la,p,m}(w,w'):=\beta_{j,p}(\la)r^{-2s-4j}
F_{j,\la}(r\omega_m){v'}^{\frac{s}{2}+j}\rho_k^N(w')\chi^k_\delta\]
is analytic for $p\geq 2N$ in $w$ in some complex
neighbourhood $U_\eps=\{ w\in \cc^{n-k+1};\textrm{dist}(w,U)\leq
\eps\}$ of $U=\textrm{supp}\nabla\chi^k_L$ and can be bounded there,
for $m\not=0$ and some $C>0$ independent of $w'$, by
\begin{equation}\label{estimthetaj}
|\theta_{j,\la,p,m}(w,w')|\leq
C^{N+p+j}\frac{e^{-\frac{\eps_0}{2}\cjg\omega_m\cjd}}{N^{2p+\Re(\la)}}
\rho_k(w')^{N+\Re(s)+2j}(N+2j)^{\Re(s)+2j}
\end{equation}
if $\Re(s)+2j\geq 0$ and by
\begin{equation}\label{estimthetaj1}
|\theta_{j,\la,p,m}(w,w')|\leq
C^{N+p+j}\frac{e^{-\frac{\eps_0}{2}\cjg\omega_m\cjd}}{N^{2p+\Re(\la)}}
\rho_k(w')^{N+\Re(s)+2j}\max\left(\cjg\omega_m\cjd^{\Re(s)+2j},
\left(\frac{\cjg\omega_m\cjd^2}{N+2j}\right)^{\Re(s)+2j}\right)
\end{equation}
if $\Re(s)+2j<0$. Indeed, we have from (\ref{rcomplex})
\[C^{-1}{v'}^\demi\rho_k(w')^{-1}\leq \Re(r^2)\leq |r|^2\leq C\Re(r^2)\leq
C{v'}^\demi\rho_k(w')^{-1}\]
\[ C^{-1}{v'}^{\frac{1}{4}}\rho_k(w')^{-\demi}\leq (\Re(r^2))^\demi\leq \Re(r)\leq
(\Re(r^2))^\demi+C{v'}^{\frac{1}{4}}\rho_k(w')^{-\demi}\leq
C{v'}^{\frac{1}{4}} \rho_k(w')^{-\demi}\] for $w\in U_\eps,
w'\in\supp(\chi^k_\delta)$, hence
\[\rho_k(w')^N|{v'}^{\frac{s}{2}+j}r^{-2s-4j}|\leq
C^{N+j}\rho_k(w')^{N+\Re(s)+2j}\]
\[\rho_k(w')^N\left|\frac{{v'}^{\frac{s}{2}+j}r^{-s-2j}}{(C\Re(r))^
{|\Re(s)+2j|}}\right| \leq
C^{N+j}\rho_k(w')^{N+\Re(s)+2j}({v'}^{-\demi}\rho_k(w'))
^{|\frac{\Re(s)}{2}+j|-(\frac{\Re(s)}{2}+j)}\] and
(\ref{estimthetaj}),(\ref{estimthetaj1}) are obtained in view of
(\ref{ks}), (\ref{estimalphaj}), the bound
${v'}^{-\demi}\rho_k(w')\leq 1$ and the uniform estimate
\[
(\Re(r)|\omega_m|)^{\Re(s)+2j}
e^{-\frac{\Re(r)|\omega_m|}{2}} \leq
(C\cjg\Re(s)+2j\cjd)^{|\Re(s)+2j|}.
\]
By Cauchy formula and
(\ref{rhomajdelta}), we deduce for $|\alpha|\leq 4N$, $w\in U$ and
$\Re(s)+2j\geq 0$
\begin{equation}\label{estimfin}
\frac{|\pl_w^\alpha\theta_{j,\la,p,m}(w,w')|}{\rho_k(w')^{n}}\leq
C^{N+j+p}\frac{\delta^{N+\Re(s)+2j}}{N^{2p-|\alpha|-2j}}
e^{-\frac{\eps_0\cjg\omega_m\cjd}{2}}
\end{equation}
whereas for $\Re(s)+2j<0$
\begin{equation}\label{estimfin2}
\frac{|\pl_w^\alpha\theta_{j,\la,p,m}(w,w')|}{\rho_k(w')^{n}}\leq
\frac{C^{N+j+p}\delta^{N+\Re(s)+2j}}{
N^{2p+\Re(\la)-|\alpha|}}e^{\frac{-\eps_0\cjg\omega_m\cjd}{2}}
\max\left(\cjg\omega_m\cjd^{\Re(s)+2j},
\left(\frac{\cjg\omega_m\cjd^2}{N+2j} \right)^{\Re(s)+2j}\right)
\end{equation}
For the case $m=0$, we obtain the same bound as
(\ref{estimfin}) by using
\[\left|\beta_{j,p}(\la)\frac{\alpha_{j,n-k}(\la-\frac{k}{2})}{\alpha_{j,n}(\la)}
\right|\leq C^{N+p+j}N^{-2p+2j}\] for $|\la|\leq \frac{N}{2}$ and
$\textrm{dist}(\la,\frac{k}{2}-\nn_0)>\frac{1}{8}$. Using
(\ref{estimprodq}), (\ref{estimprodq2}) with $M=N$ and $p\leq 2N$,
(\ref{estimfin}), (\ref{estimfin2}), (\ref{quasi}) and again the
bound
\[\cjg\omega_m\cjd^{\Lambda}e^{-\frac{\eps_0\cjg\omega_m\cjd}{4}} \leq
(C\Lambda)^\Lambda\]
for all $\Lambda>0$, we can conclude that
\[\left|\left|(\Delta_{\Omega}+1)^N\rho_k^{-N}\phi_LK_{2N,m}(\la)\rho_k^N\right|\right|
_{\mc{L}(\mc{H}_m,L^2(\Omega))} \leq
\frac{\delta^{\Re(s)+N}(CN)^{2N}} {e^{\frac{\eps_0\cjg\omega_m\cjd}{4}}}
\max\Big(1,(\cjg\omega_m\cjd N^{-1})^{\Re(s)-|\Re(s)|}\Big)\] and the
same estimate for $L_{2N,m}(\la)$. We just recall that the
singular values of $(1+\Delta_\Omega)^{-N}$ on $L^2(\Omega)$ satisfy
\[\mu_l((1+\Delta_\Omega)^{-N})\leq
(Cl)^{\frac{2N}{n-k+1}}\] and that
\begin{equation}\label{prodvs}
\mu_l(AB)\leq \mu_l(A)||B||
\end{equation}
if $A$ is trace class and $B$ bounded to show that
\[\mu_l\Big(\phi_LK_{2N,m}(\la)+L_{2N,m}(\la)\Big)\leq e^{-\frac{\eps_0\cjg\omega_m\cjd}{4}}
\delta^{\Re(s)+N}\left(Cl^{-\frac{1}{n-k+1}}N\right)^{2N}
\max\Big(1,\Big(\frac{\cjg\omega_m\cjd}{N}\Big)^{\Re(s)-|\Re(s)|}\Big)\]
on $\rho_k^N\mc{H}_m$. In view of (\ref{comprhok}), (\ref{comprho}) and
(\ref{prodvs})
this gives
\[\mu_l\Big(\phi_LK_{2N,m}(\la)+L_{2N,m}(\la)\Big)\leq e^{-\frac{\eps_0\cjg\omega_m\cjd}{4}}
\delta^{\Re(s)}\left(Cl^{-\frac{1}{n-k+1}}N\right)^{2N}
\max\Big(1,\Big(\frac{\cjg\omega_m\cjd}{N}\Big)^{\Re(s)-|\Re(s)|}\Big)\]
on $\rho_\delta^N\mc{H}_m$.

By taking $M=0$ in (\ref{estimprodq}), the previous estimates also
show that for $\Re(\la)>\ndemi$
\[\left|\left|\phi_LK_{2N,m}(\la)+L_{2N,m}(\la)\right|\right|_{\mc{L}(\rho^N_\delta\mc{H}_m)}\leq
C^N\delta^{\Re(s)}e^{-\frac{\eps_0\cjg\omega_m\cjd}{4}}\]
and (\ref{majklnorm}) is then easily deduced.
\qed\\

In a second step, we are going to control the singular values of
the terms $K^\sharp_{2N,m}(\la)$ and
$L^\sharp_m(\la)$ on $\rho_\delta^N\mc{H}_m$.

\begin{lem}\label{singvalkdiese}
$K^\sharp_{2N,m}(\la)$ and
$L^\sharp_m(\la)$ are trace class on $\rho_\delta^N\mc{H}_m$ if
$\delta>0$ is chosen small enough and $|\la|\leq \frac{N}{2}$,
$\textrm{dist}(\la,\frac{k}{2}-\nn_0)>\frac{1}{8}$. Moreover, their singular
values satisfy
\[\mu_l\Big(K^\sharp_{2N,m}(\la)+L_m^\sharp(\la)\Big)\leq
e^{-\frac{\eps_0\cjg\omega_m\cjd}{4}}
\delta^{\Re(s)}\left(Cl^{-\frac{1}{n-k+1}}N\right)^{2N}
\max\Big(1,\Big(\frac{\cjg\omega_m\cjd}{N}\Big)^{\Re(s)-|\Re(s)|}\Big)\]
there for some $\eps_0>0$
and if $\la_N=\frac{N}{4}$ we have 
\begin{equation}\label{majkldiesenorm}
\left|\left|K^\sharp_{2N}(\la_N)+L^\sharp(\la_N)
\right|\right|_{\mc{L}(\rho^N_\delta L^2(M_k))}\leq (C\delta)^{\frac{N}{4}}
\end{equation}
\end{lem}
\textsl{Proof}: we recall that
\[\rho_k^N(w)K^\sharp_{2N,m}(\la;w,w')\rho_k^N(w')=
\rho_k^Nv^{\frac{s}{2}+2N+1}\phi_L[\Delta_{X_k},\psi_L^k]
\sum_{j\geq 2N+1}v^{j-2N-1}\theta_{j,\la,j,m}(w,w')\] provided the
sum converges. Taking advantage of the estimates (\ref{estimfin}),
we find for $w\in U_\eps$, $|\la|\leq \frac{N}{2}$, $\Re(s)+2j\geq 0$
and $p\geq 0$
\begin{equation}\label{vpthet}
|v^{p}\theta_{j,\la,j,m}(w,w')|\leq C^{N+p}\delta^{\Re(s)+N}
e^{-\frac{\eps_0\cjg\omega_m\cjd}{2}}(C\delta)^j\rho_k(w')^{n}
\end{equation}
which proves that the sum converges if $\delta$ is chosen small
enough and we obtain
\[|\rho_k^N(w)K^\sharp_{2N,m}(\la;w,w')\rho_k^N(w')|\leq
(C\delta)^N\delta^{\Re(s)+N}e^{-\frac{\eps_0\cjg\omega_m\cjd}{4}}\rho_k(w')^{n}.\]
Using the arguments of Lemma \ref{singvalk}, it is straightforward
to check that
\[\mu_l(K^\sharp_{2N,m}(\la))\leq
(C\delta)^N e^{-\frac{\eps_0\cjg\omega_m\cjd}{4}}l^{-2}.\]

For $L^\sharp(\la)$ the method is similar, we recall that
\[\rho_k^{-N}(w)L^\sharp_{m}(\la;w,w')\rho_k^N(w')=
\rho_k^{-N}\psi_L^k[\Delta_{X_k},\phi_L]v^{\frac{s}{2}}
\sum_{j=0}^\infty v^{j}\theta_{j,\la,j,m}(w,w').\] Using
(\ref{vpthet}) for $\Re(s)+2j\geq 0$ and (\ref{estimfin2}) for
$\Re(s)+2j\leq 0$, we find that
\[|\pl^N_w(v^j\theta_{j,\la,j,m}(w,w'))|\leq
e^{-\frac{\eps_0\cjg\omega_m\cjd}{2}}
(CN)^{2N}\delta^{\Re(s)+N}(C\delta)^j\max(1,(\cjg\omega_m\cjd N^{-1})
^{\Re(s)-|\Re(s)|})\rho_k(w')^n\]
for $w$ in the complex neighbourhood $U_\eps$
of $\supp(\nabla\chi_L^k)$. One deduces that the sum $L^\sharp_{m}(\la)$
converges for small $\delta$ and the arguments of Lemma
\ref{singvalk} yield the bound
\[\mu_l(L^\sharp_m(\la))\leq
e^{-\frac{\eps_0\cjg\omega_m\cjd}{4}}
\delta^{\Re(s)}\left(Cl^{-\frac{1}{n-k+1}}N\right)^{2N}
\max\Big(1,(N^{-1}\cjg\omega_m\cjd)^{\Re(s)-|\Re(s)|}\Big)\]
on $\rho^N_\delta\mc{H}_m$.
These estimates on the singular values also imply
(\ref{majkldiesenorm}).
\qed\\

Lemmas \ref{singvalk} and \ref{singvalkdiese} clearly prove that
\[\sum_{m\in\zz^k}\sum_{l\in\nn}\mu_l(\mc{K}^k_{N,m}(\la))<\infty\]
and we have
\[\prod_{m\in \zz^k}\prod_{l\in\nn}(1+q\mu_l(\mc{K}^k_{N,m}(\la)))\leq
\exp\left(\sum_{m\in\zz^k}\sum_{l\in\nn}\log\left(1+(\Lambda l^{-\frac{1}{n-k+1}})^{2N}\right)\right).\]
with $\Lambda:=Cq\delta^{-1}e^{-\frac{\eps_0\cjg\omega_m\cjd}{8N}}N
\max(1,(N\cjg\omega_m\cjd^{-1})^\demi)$. Now, we use
\begin{eqnarray*}
\sum_l\log\left(1+(\Lambda l^{-\frac{1}{n-k+1}})^{2N}\right)&\leq
&\int_0^\infty \log\left(1+(\Lambda t^{-\frac{1}{n-k+1}})^{2N}\right)dt\\
& \leq &\Lambda^{n-k+1}\int_0^\infty \log(1+t^{-\frac{2N}{n-k+1}})dt\\
&\leq & CN\Lambda^{n-k+1}\int_1^\infty t^{-\frac{3}{2}}dt\\
& \leq & C_{q,\delta} e^{-\frac{\eps_0\cjg\omega_m\cjd}{8N}}
N^{n-k+2}\max(1,(N\cjg\omega_m\cjd^{-1})^\demi).
\end{eqnarray*}
for some $C_{q,\delta}>0$. 
Finally, since $\cjg\omega_m\cjd\geq C|m|$ for some $C>0$ depending on
$A_k\in GL_k(\rr)$ we have
\[\sum_{\substack{m\in\zz^k\\ |\omega_m|\geq N}}
e^{-\frac{\eps_0\cjg\omega_m\cjd}{8N}}  \leq \sum_{j\in\nn_0}\sum_{j
\leq |m|\leq j+1} e^{-\frac{j}{CN}} \leq C\int_0^\infty
t^{k-1}e^{-\frac{t}{CN}}dt \leq CN^{k}\] and
\[\sum_{\substack{m\in\zz^k\\ \cjg\omega_m\cjd\leq N}}
e^{-\frac{\eps_0\cjg\omega_m\cjd}{8N}}\sqrt{\frac{N}{\cjg\omega_m\cjd}} \leq
\sum_{j\in\nn_0}\sum_{j \leq |m|\leq j+1} e^{-\frac{j}{CN}}
\sqrt{\frac{N}{\cjg j\cjd}} \leq C\int_0^\infty
t^{k-1}\sqrt{\frac{N}{t}} e^{-\frac{t}{CN}}dt \leq CN^{k}.\]
This proves that we can find $C_{\delta,q}>0$  such that
\[\prod_{m\in\zz^k}\prod_{l\in\nn}(1+q\mu_l(\mc{K}^k_{N,m}(\la)))
\leq C_{\delta,q} e^{C_{\delta,q} N^{n+2}}\] for $|\la|\leq \frac{N}{2}$,
$\textrm{dist}(\la,\frac{k}{2}-\nn_0)>\frac{1}{8}$, thus (\ref{detkcusp})
is obtained.\\

The bound (\ref{normlan}) is a consequence of Lemmas \ref{singvalk} and
\ref{singvalkdiese}.\\

To conclude the proof of Theorem \ref{parametrixkcusp}, it remains
to prove the  
\begin{lem}\label{contek}
The operator $\mc{E}^k_N(\la)$ is continuous from
$\rho_\delta^NL^2(M_k)$ to $\rho_\delta^{-N}L^2(M_k)$. 
\end{lem}
\textsl{Proof}: except $\rho_\delta^NE_0(\la)\rho^N_\delta$,
we have seen that the other terms in the expression of $\rho_\delta^N\mc{E}^k_N(\la)\rho^N_\delta$
have Schwartz kernels in $L^2(M_k\x M_k)$ and thus are bounded
on $L^2(M_k)$.
To deal with $\rho_\delta^NE_0(\la)\rho^N_\delta$, we take $J=2N$ in (\ref{rkm}) 
and first show that
\begin{equation}\label{kernel}
(\rho_\delta(w)\rho_\delta(w'))^Nd^{s+2N}\int_{\rr^k}
e^{-ir\omega_m.z}(1+|z|^2)^{-s-2N}G_{2N,n}(s,d(1+|z|^2)^{-1})dz
\end{equation}
are the kernels of bounded operators on $\mc{H}_m$ with norms
uniformly bounded with respect to $m$ when $|s|\leq \frac{N}{2}$.
We know that $G_{2N,n}(s,\tau)$ is smooth for $\tau\in[0,\demi)$
thus (\ref{kernel}) is square integrable in
$(M_k\x M_k)\setminus \{d>\frac{1}{8}\}$ with norm bounded
by $C_N$ for $|s|\leq \frac{N}{2}$. Now using  \cite[Prop. A.1]{P}, we deduce that
(\ref{kernel}) is bounded for $\{d>\frac{1}{8}\}$ by
$\varphi(d_{\hh^{n-k+1}}(w,w'))$ for a function $\varphi>0$
satisfying
\[\int_0^1\varphi(\tau)\tau^{n-k}d\tau\leq C_N\]
for $|s|\leq \frac{N}{2}$.
Therefore, Proposition B.1 of \cite{P} allows to conclude that (\ref{kernel}) is bounded
on $\mc{H}_m$ uniformly with respect to $m$.
Now it just remains to show the boundedness of the operator on $\mc{H}_m$ whose kernel is 
$\rho_k(w)^N\rho_k(w')^Nd^{s+2j}F_{j,\la}(r|\omega_m|)$ with a uniform bound with respect to $m$. 
First we use that 
\[\rho_k(w)^N\rho_k(w')^Nd^{\Re(s)+2j}\leq C_N\rho_k(w)^{\frac{N}{2}}\rho_k(w')^{\frac{N}{2}}\]
if $|s|\leq \frac{N}{2}$, which is straightforward for $\Re(s)+2j>0$ and comes easily from the estimate
\[d^{-1}<C\left(\frac{1}{\rho_k(w)x'}+\frac{1}{\rho_k(w')x}\right)\]
if $\Re(s)+2j<0$. Now to bound $F_{j,\la}(r|\omega_m|)$ we see that if $r|\omega_m|>1$
then 
\[(r|\omega_m|)^{s+2j}K_{-s-2j}(r|\omega_m|)|\leq C_{N}\]
in view of the bound (\ref{ks}), whereas if $r|\omega_m|<1$ we can use that
$K_s(z)=z^s\phi_{s}(z^2)+z^{-s}\phi_{-s}(z^2)$ for some function $\phi_s(z^2)$
smooth on $z\in[0,1]$ and observe that
\[|F_{j,\la}(r|\omega_m|)|\leq|\phi_{-s-2j}(r|\omega_m|)|+|(r|\omega_m|)^{2s+4j}\phi_{s+2j}(r|\omega_m|)|\]
is bounded by $C_N$ if $\Re(s)+4j>0$ and by $C_Nr^{2\Re(s)+4j}\leq C_N
d^{-\Re(s)-2j}(\rho_k(w)\rho_k(w'))^{-\frac{N}{2}}$ if $\Re(s)+2j<0$. 
This proves that in all cases we have
\[|\rho_k(w)^N\rho_k(w')^Nd^{s+2j}F_{j,\la}(r|\omega_m|)|\leq 
C_N(\rho_k(w)\rho_k(w'))^{\frac{N}{2}}\] 
and  we conclude that $\rho_k(w)^N\rho_k(w')^Nd^{s+2j}F_{j,\la}(r|\omega_m|)$ is the kernel of a Hilbert Schmidt
operator on $\mc{H}_m$ with norm uniformly bounded with respect to $m$, this
achieves the proof of the Lemma.
\qed\\

The Proposition \ref{parametrixkcusp} is then proved.
\qed\\

\subsection{The maximal rank cusps}
We recall that $X_n=\Gamma_n\backslash \hh^{n+1}$ is a quotient by
a group of translations acting on $\rr^n$. The lattice of translations
$\Gamma_n$ acting on $\rr^n$ is the image of
the lattice $\zz^k$ by a map $A_n\in GL_n(\rr)$.
By using a Fourier decomposition on the torus $T^n=\Gamma_n\backslash\rr^n$ and
conjugating by $x^{\frac{n}{2}}$,
the operator $\Delta_{X_n}-\la(n-\la)$ acts on
\[L^2(X_n)=\bigoplus_{m\in\zz^n}\mc{H}_m,\quad
\mc{H}_m\simeq L^2(\rr^+,x^{-1}dx)\]
as a family of operators
\[P_m(\la):=-(x\pl_x)^2+x^2|\omega_m|^2+s^2\]
where $\omega_m=2\pi \trans (A_n^{-1})m$ for $m\in\zz^n$
and $s:=\la-\frac{n}{2}$ the shifted spectral parameter.
By elementary Sturm-Liouville theory (see \cite[Lem. 3.1]{G}),
we find that the resolvent $R_{X_n}(\la)=(\Delta_{X_n}-\la(n-\la))^{-1}$
for the Laplacian on $X_n$ is for $\Re(\la)>\ndemi$
\begin{equation}\label{resolventxn}
R_{X_n}(\la)=\bigoplus_{m\in\zz^n}R_m(\la) \textrm{ on }
L^2(X_n)=\bigoplus_{m\in\zz^n}\mc{H}_m
\end{equation}
with
\[R_m(\la;x,x')=-K_{-s}(|\omega_m|x)I_s(|\omega_m|x')H(x-x')-
K_{-s}(|\omega_m|x')I_s(|\omega_m|x)H(x'-x), \quad m\not=0
\]
\[R_0(\la;x,x')=(2s)^{-1}e^{-s|\log(x/x')|}\]
where $H$ is the Heaviside function, $K_s$ is defined in
(\ref{besselk}) and $I_s$ is the modified Bessel function.

Now we construct a parametrix for $\Delta_X-\la(n-\la)$ on
the end $I_n^{-1}(M_n)$ of our manifold $X$. Notice however that
better estimates could be obtained for this part (see e.g. \cite{GZ1})
since the problem is essentially reduced to the one-dimensional case.

\begin{prop}\label{parametrixncusp}
There exist some bounded operators
\[\mc{E}_N^n(\la): \rho_\delta^NL^2(M_n)\to \rho_\delta^{-N}L^2(M_n)\]
\[\mc{K}_N^n(\la): \rho_\delta^NL^2(M_n)\to \rho_\delta^NL^2(M_n)\]
holomorphic in $\{\Re(\la)>-\frac{N+1}{2},\la\not=\ndemi\}$ with at most a simple pole
at $\ndemi$ and such that
\[(\Delta_{M_n}-\la(n-\la))\mc{E}_N^n(\la)=\chi_\delta^n+\mc{K}_N^n(\la),\]
$\mc{K}^n_N(\la)$ is trace class on $\rho_\delta^NL^2(M_n)$ and for $q>0$ 
there exists $C_{\delta,q}>0$ such that for $|\la|\leq \frac{N}{2}$
\[\det(1+q|\mc{K}^n_N(\la)|)\leq e^{C_{\delta,q}\cjg N\cjd^{n+2}},\]
the determinant being on $\rho_\delta^NL^2(M_n)$.
Moreover, for $\la_N=\frac{N}{4}$,
\[
||\mc{K}^n_N(\la_N)||_{\mc{L}(\rho_\delta^NL^2(M_n))}\leq
(C\delta)^{\frac{N}{4}}.
\]
\end{prop}
\textsl{Proof}: let us set
\[\mc{E}_N^n(\la):=\chi^n_LR_{X_n}(\la)\chi_\delta^n,
\quad \mc{K}_N^n(\la):=[\Delta_{X_n},\chi^n_L]R_{X_n}(\la)\chi_\delta^n\]
and check that this choice satisfies the announced properties.
The boundedness of $\mc{E}_N^n(\la)$ is obtained by Schur's lemma.
To show that $\mc{K}_N^n(\la)$ is trace class on $\rho_\delta^NL^2(M_n)$
and to estimate its singular values,
we analyze each $[x\pl_x,\phi_L]R_m(\la)\chi_\delta^n$
on $\rho_\delta^N\mc{H}_m$ and
use the same arguments as for the non-maximal rank cusps.
Since we do not need the optimal estimates for the singular values,
it suffices to control the derivatives of $R_m(\la;x,x')$
outside the diagonal.
This is easily obtained from the formulae
\[
I_k(z)=\frac{1}{\pi}
\int_0^\pi e^{z\cos(u)}\cos(ku)du-\frac{\sin(k\pi)}{\pi}\int_0^\infty
e^{-z\cosh(u)-ku}du
\]
\[
K_{-k}(z)=\int_0^\infty\cosh(ku)e^{-z\cosh(u)}du.
\]
the analyticity in $z$ and Cauchy's formula as for the non-maximal rank cusps.
Then a straightforward calculus shows that the singular values of
$[x\pl_x,\phi_L]R_m(\la)\chi_\delta^n$ on $\rho_\delta^N\mc{H}_m$ satisfy
\[\mu_l([x\pl_x,\phi_L]R_m(\la)\chi_\delta^n)\leq
e^{-\frac{\eps_0\cjg\omega_m\cjd}{4}}
\delta^{\Re(s)}\left(Cl^{-1}N\right)^{2N}
\max\Big(1,(N^{-1}\cjg\omega_m\cjd)^{\Re(s)-|\Re(s)|}\Big)\]
for some $\eps_0>0$ and the arguments of Proposition \ref{parametrixkcusp} allow
to complete the proof.
\qed\\

\subsection{The regular neighbourhoods}
For this part of the parametrix, we use the work
of Guillop\'e-Zworski \cite{GZ2} and deduce the following
\begin{prop}\label{parametrixreg}
There exists some bounded operators
\[\mc{E}_N^r(\la): \rho_\delta^NL^2(M_r)\to \rho_\delta^{-N}L^2(M_r)\]
\[\mc{K}_N^r(\la): \rho_\delta^NL^2(M_r)\to \rho_\delta^NL^2(M_r)\]
meromorphic in $\Re(\la)>-\frac{N+1}{2}$ with simple poles
at $-j$ (with $j\in\nn_0$) of ranks uniformly
bounded by $C(j+1)^{n+1}$ such that
\[(\Delta_{M_r}-\la(n-\la))\mc{E}_N^r(\la)=\chi_\delta^r+\mc{K}_N^r(\la),\]
$\mc{K}^r_N(\la)$ is trace class on $\rho_\delta^NL^2(M_r)$ and for
$q>0$ and $|\la|\leq \frac{N}{2}$ and $\textrm{dist}(\la,-\nn_0)>\frac{1}{8}$ 
there exists $C_{\delta,q}>0$ such that
\[
\det(1+q|\mc{K}^r_N(\la)|)\leq e^{C_{\delta,q}\cjg N\cjd^{n+2}}, \]
the determinant being on $\rho_\delta^NL^2(M_r)$.
Moreover, for $\la_N=\frac{N}{4}$,
\[
||\mc{K}^r_N(\la_N)||_{\mc{L}(\rho_\delta^NL^2(M_r))}\leq
(C\delta)^{\frac{N}{4}}.
\]
\end{prop}
\textsl{Proof}:  we begin by defining a new cut-off function 
$\chi^r_{L,\delta}(x,y)=\phi_L(Cx\delta^{-1})\psi^r_L(y)$ with 
$C$ chosen so that $\chi^r_{L,\delta}=1$ on $\supp(\chi^r_L)$.
It suffices now to use the construction of \cite[Prop. 3.1]{GZ2} with the cut-off functions
$\chi_{L,\delta}^r$ and $\chi_\delta^r$
(i.e. the functions $\chi_1^\delta, \chi_2^\delta$ of \cite{GZ2} 
are $\chi_{L,\delta}^r,\chi_\delta^r$ here).
Note that in their construction, Guillop\'e and Zworski used for
$\chi_2^\delta$ a function having a product structure like
$\chi_{L,\delta}^r$ but it is not difficult to see that 
our $\chi^r_\delta$ suits as well by using that
$C^{-1}x\leq \rho_\delta\leq Cx\delta^{-1}$ in $M_r$ for some $C>0$.
The end of the proof is given by Proposition 4.1 of \cite{GZ2}.
\qed\\

\section{Bounds on resonances}

Combining the Propositions
\ref{parametrixkcusp}-\ref{parametrixncusp}-\ref{parametrixreg}
and the same kind of arguments used by Guillop\'e-Zworski
\cite{GZ2}, we can prove the Theorems.\\

\textsl{Proof of Theorem \ref{prolongmer}}: the first thing is to
construct the final parametrix. We define for $\Re(\la)>-\frac{N+1}{2}$
\[\mc{E}_N(\la):= \chi^i_{L,\delta}R(\la_N)\chi^i_\delta+\sum_{k=1}^n
(I_k)^*\mc{E}^k_N(\la)(I_k)_*+(I_r)^*\mc{E}^r_N(\la)(I_r)_*\]
\[\mc{K}_N(\la):=\mc{K}_N^i(\la)+
\sum_{k=1}^n(I_k)^*\mc{K}^k_N(\la)(I_k)_*+(I_r)^*\mc{K}^r_N(\la)(I_r)_*\]
\[\mc{K}^i_N(\la):=[\Delta_X,\chi^i_{L,\delta}]R(\la_N)\chi^i_\delta+
(\la(n-\la)-\la_N(n-\la_N))\chi^i_{L,\delta}R(\la_N)\chi^i_\delta\]
and we get by construction
\[(\Delta_X-\la(n-\la))\mc{E}_N(\la)=1+\mc{K}_N(\la).\]
Moreover we deduce from the Propositions
\ref{parametrixkcusp}, \ref{parametrixncusp} and \ref{parametrixreg}
that $\mc{K}_N(\la)^{n+2}$ is trace class on $\rho_\delta^NL^2(X)$
such that $||\mc{K}_N(\la_N)||_{\mc{L}(\rho_\delta^NL^2(X))}\leq \demi$ if
$\delta$ is chosen small and $N$ large, and
$\mc{E}_N(\la)$ is bounded from $\rho_\delta^NL^2(X)$ to
$\rho_\delta^{-N}L^2(X)$.
Consequently, $(1+\mc{K}_N(\la_N))$ is invertible on
$\rho_\delta^NL^2(X)$ and Fredholm analytic theory allows to invert
$(1+\mc{K}_N(\la))$ meromorphically with finite rank poles on the same
Hilbert space, which gives the analytic continuation of $R(\la)$ to $\{\Re(\la)>-\frac{N+1}{2}\}$ 
as a family of operators in $\mc{L}(\rho^NL^2(X),\rho^{-N}L^2(X))$, thus to $\cc$ 
from $L^2_{comp}(X)$ to $L^2_{loc}(X)$ since $N$ can be chosen arbitrarily large.
\qed\\

\textsl{Proof of Theorem \ref{nbres}}:
we define the determinant
\[D_N(\la):=\det(1+\mc{K}_N(\la)^{n+2})\]
on $\rho_\delta^NL^2(X)$. This is a meromorphic function in
$\Re(\la)>\frac{N+1}{2}$ such that
the resonances of $\Delta_X$ are contained in the set of
zeros of $D_N$ with multiplicities and $\demi(n-\nn)$
with multiplicity of $\demi(n-j)$ bounded by $v_j(D_N)+C(j+1)^{n+1}$
where $v_j(D_N)$ is the order of the zero (or pole) $j$ for $D_N$
(see the appendix of \cite{GZ2} for the multiplicity). Moreover
by fixing $\delta$ small enough, it is clear that  
$|\det(1+\mc{K}_N(\la_N)^{n+2})|>\demi$ for $N$ large.\\

$\delta$ is now fixed as before, we then use Lemma 6.1 of \cite{GZ1} and deduce that
\[|D_N(\la)|\leq q\sum_{k=1,\dots,n,r}
\det(1+q|\mc{K}^k_N(\la)|)^{q}+
q\det(1+q|\mc{K}_N^i(\la)|^{n+2})^{q}\]
for some $q>0$ independant of $N$.
Moreover, it is straightforward to see that there exists $C>0$ such that
\[\mu_l(\mc{K}^i_N(\la))\leq
C(|\la-\la_N|+1)l^{-\frac{1}{n+1}} \]
which combined with the Propositions
\ref{parametrixkcusp}, \ref{parametrixncusp} and \ref{parametrixreg} shows
that
\[|D_N(\la)|\leq Ce^{CN^{n+2}}\]
for $|\la|\leq \frac{N}{2}$ and $\textrm{dist}(\la,\demi(n-\nn))>\frac{1}{8}$.
To complete the proof, it suffices to multiply $D_N(\la)$ by the function
\[g_P(\la):=\la^{P2^{n+1}}\prod_{\omega\in U_{2n+4}.\demi(n-\nn)}
\left(E\left(\frac{\la}{\omega},n+2\right)\right)^{P(|2\omega|)^{n+1}}\]
as defined in Section 5 of \cite{GZ2} ($U_m$ being the set of
m-th root of the unity and $E(z,p):=(1-z)\exp(z+\dots+p^{-1}z^p)$ are 
the elementary Weierstrass functions)
and use the Lemmas 5.1 and 5.3 of this article to prove that for $P$ big enough (independant of $N$),
$D_N(\la)g_P(\la)$  is a holomorphic function bounded
by $Ce^{CN^{n+2}}$ in the disc $|\la|\leq \frac{N}{2}$ (note that 
the maximum principle is used to control the norms near the points
$\frac{n-j}{2}$ for $j\in\nn$).
In view of the discussion about the relation between resonance multiplicity
and the valuation of determinant, this completes the proof of the
Theorem by applying Jensen's lemma to $g_pD_N$
in the disc centered in $\la_N=\frac{N}{4}$
with radius $\frac{N}{2}$.
\qed\\

\end{document}